\documentclass[journal]{IEEEtran}

\usepackage{graphicx}      
\usepackage[applemac]{inputenc} 
\usepackage{amsmath}
\usepackage{amssymb}
\usepackage{amsfonts}

\newcommand{\bma}{\left(\begin{array}}
\newcommand{\ema}{\end{array}\right)}

\newcommand{\bsubeq}{\begin{subequations}}
\newcommand{\esubeq}{\end{subequations}}

\newcommand{\biz}{\begin{itemize}}
\newcommand{\eiz}{\end{itemize}}

\newcommand{\benu}{\begin{enumerate}}
\newcommand{\eenu}{\end{enumerate}}

\newcommand{\ece}{\end{center}}
\newcommand{\bce}{\begin{center}}
\newcommand{\bem}{\begin{em}}
\newcommand{\eem}{\end{em}}

\newcommand{\bpm}{\begin{pmatrix}}
\newcommand{\epm}{\end{pmatrix}}

\newtheorem{thm}{{\bf Theorem}}

\usepackage{color}

\usepackage{bm}

\usepackage{amssymb}
\usepackage{amsfonts}
\usepackage{amsmath}
\usepackage{graphicx}
\usepackage{epstopdf}
\usepackage{color}

\usepackage{bbm}

\newtheorem{theorem}{Theorem}

\newtheorem{lemma}[theorem]{Lemma}

\newdimen\dummy

\allowdisplaybreaks

\begin{document}
%
\title{State Feedback Stabilization of the Linearized Bilayer  \textit{Saint-Venant}  Model}
%
%
%

\author{Ababacar Diagne, Shuxia Tang, Mamadou Diagne, and Miroslav Krstic
\thanks{A. Diagne is  with the Division of Scientific Computing, Department of Information Technology,  Uppsala University,
Box 337, 75105 Uppsala, Sweden (email: ababacar.diagne@it.uu.se).

S. Tang (corresponding author) and M. Krstic are with the Department of Mechanical and Aerospace Engineering, University of California, La Jolla, CA 92093-0411. 
(e-mail: sht015@ucsd.edu; krstic@ucsd.edu)

 M. Diagne is  with the Department of Mechanical Engineering, University of Michigan G.G. Brown Laboratory
2350 Hayward
Ann Arbor MI 48109, USA (e-mail: mdiagne@umich.edu)
}
}

\maketitle

\begin{abstract}
	We consider the problem of stabilizing the  bilayer  \textit{Saint-Venant}  model, which is a coupled system of two rightward
and two
leftward convecting transport partial 
differential
equations (PDEs).
In the stability proofs, we employ a  Lyapunov
function in which the parameters need   to be successively determined. To the best of the authors' knowledge, 
this is the first time this kind of Lyapunov function is employed, and this result is the first one on  the stabilization of the linearized bilayer  \textit{Saint-Venant}
model. Numerical simulations of the bilayer  \textit{Saint-Venant}  problem are also provided to verify the result.
\end{abstract}


%
\IEEEpeerreviewmaketitle

%
%
%
%
\section{Introduction}
Population and economic growth are changing  social values about the importance
of water and the expansion of the energy sector  will continue to drive growing  demand for water resources. These trends are placing greater pressure on existing water allocations,
heightening the importance of water management and conservation for the
sustainability of irrigated agriculture. During the past decades,  several efforts have been made by engineers and  researchers towards the design of control methodologies for the real-time monitoring  of
irrigation canals. The global challenge motivating these studies is becoming increasingly
obvious  due to the less than desirable performance of both manually
automatically  controlled water management infrastructures.

In this paper,  the $1$D two-layer Saint-Venant model that consists of the superposition
of two immiscible  fluids  with different constant densities is presented. The derivation of the bilayer  model can be found in \cite{Audusse2011} and
 \cite{bouchut2008}. The later one develops a stable well-balanced time-splitting scheme for a type of bilayer Saint-Venant model which  satisfies a fully discrete entropy inequality.
One can find  some results on mathematical analysis of  the related  problem in  \cite{Narbona2011} and  \cite{Munoz2003}. To the best of the authors' knowledge,  the relevant   control related problem for this  application has not been investigated in the existing contributions. 

PDE backstepping control approach  has been  successfully
 employed for the feedback stabilization of various  classes  of PDEs \cite{Dimeglio2013, krstic2008boundary}.  In the present work,  a general system, which  consists of $m$ rightward
and $n$ leftward
transport PDEs with spatially varying coefficients, is exponentially stabilized by $m$ boundary input backstepping controllers.
Our  backstepping controller design idea can be referred to  \cite{LongHu2015}, in which the stabilization problem 
for the general coupled heterodirectional  system of hyperbolic types with an arbritrary number of equations convecting  in both directions is definitely solved.  In our stability proofs, we employ a Lyapunov
function in which the parameters need   to be successively determined. 
Then, applying this general stabilization result to the 1D bilayer  \textit{Saint-Venant} problem, which
 consists of two rightward and two leftward convecting transport PDEs $(n=m=2)$, we  achieve exponential stabilization 
 with two boundary input controllers.  

 This paper is organized as follows. In  Section \ref{sec1}, we state the control problem. The
 1D bilayer  \textit{Saint-Venant} model
 is
 first formulated based on its physical
 description, and then a linearized version around a steady state
 is presented.
 Section 3 is dedicated to the state feedback backstepping
 controller design of a more general system, for which and exponential stability is achieved for the closed-loop control system.  This result could serve as a full theoretical result by itself, and it can be immediately utilized  for the linearized bilayer  \textit{Saint-Venant} model. Numerical simulations are
provided in Section 4. Finally in Section 5, a conclusion is
 presented and some perspectives are discussed.

\section{Problem Statement}\label{sec1}
\subsection{The 1D nonlinear bilayer  \textit{Saint-Venant }model}
We consider  the stabilization problem of a 1D   two-layer  \textit{Saint-Venant} model, which governs the dynamic of two superposed immiscible layers  of shallow
water fluids.
\begin{figure}[b]
        \centering
        \includegraphics[width=0.35\textwidth]{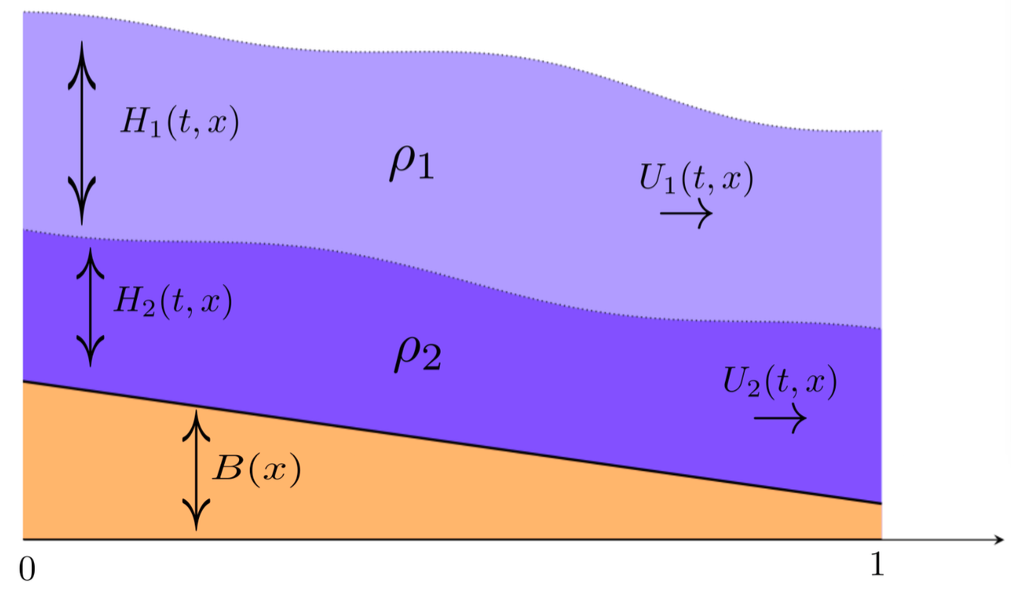}
        \label{fig:domainebicouche}
        \caption{Two-layer shallow water flows with variable topography. $H_{\rm i}$ and $U_{\rm i}$
        deote respectively the thickness of the layer $i$ and its velocity.}
\end{figure}
%
\begin{align}
\hspace{-0.1 in}\left\{
\begin{array}{ll}
\frac{\partial H_{1}}{\partial t} +  \frac{\partial (H_{1}U_{1})}{\partial
x}=0,\;\;\quad \\\\
\frac{\partial U_{1}}{\partial t} + U_{1}\frac{\partial U_{1}}{\partial x}
+ g\frac{\partial H_{1}}{\partial x} +g\frac{\partial H_{2}}{\partial x}
+ g\frac{\partial B}{\partial x}=S_{1}^{f}, 
\;\;\quad\\\\
\frac{\partial H_{2}}{\partial t} +  \frac{\partial (H_{2}U_{2})}{\partial
x}=0,\;\;\quad \\\\
\frac{\partial U_{2}}{\partial t} + U_{2}\frac{\partial U_{2}}{\partial x}
+ g\frac{\partial H_{2}}{\partial x} +g\frac{\rho_{1}}{\rho_{2}}\frac{\partial
H_{1}}{\partial x}
+ g\frac{\partial B}{\partial x}=S_{2}^{f}.
\end{array}\right.
\label{model30}
\end{align}
 In these equations,  the index $1$  refers to the upper layer and the index
$2$ to the lower one, as depicted 
in  Figure \ref{fig:domainebicouche}. 
The unknown state variables $H_i$, $U_i$ and $B$ represent respectively  the thickness of
the $i$-th layer, 
the velocity and the height of the sediment layer.
Each layer is  assumed to have a constant density  $\rho_i$, $i=1,2$ $(\rho_1<\rho_2).$
 The system contains the source terms due to the bottom topography  and the friction term.
The quantities $S_{1}^{f}$ and $S_{2}^{f}$ stand as the friction between
the two layers, and they are given by
\begin{align}
S_{1}^{f}=-C_f|U_1-U_2|(U_1-U_2)
\end{align}
and
\begin{align}
 S_{2}^{f}=rC_f|U_1-U_2|(U_1-U_2).
\end{align}

Define a vector  $ W=\begin{bmatrix} H_1, & U_1, & H_2, & U_2\end{bmatrix}^{T}$,
a ratio $ r=\frac{\rho_1}{\rho_2}$  and a map
\begin{align}
F(W)=\begin{bmatrix}
H_{1}U_{1} \\
\frac{U^{2}_{1}}{2}+g(H_1+H_2) \\
H_{2}U_{2}\\
\frac{U^{2}_{2}}{2}+g(H_2+rH_1) \\
\end{bmatrix}
\label{nota3}
\end{align}
then we could recast equation  (\ref{model30}) under the form of
\begin{align}
\frac{\partial {W}}{\partial t} +\frac{\partial F(W)}{\partial x}= S(x,W),
\label{cpmmodel3}       
\end{align}
where
\begin{align}
S(x,W)\,=\,
\begin{pmatrix}
0 & \
S_{1}^{f}-  g\frac{\partial B}{\partial x} & \
0 & \
S_{2}^{f}-  g\frac{\partial B}{\partial x} 
\end{pmatrix}^T.
\label{source1}
\end{align}
By considering the Jacobian matrix $A$ from (\ref{cpmmodel3}), we could rewrite
the equation (\ref{cpmmodel3})
into a quasilinear form as
\begin{align}
\bm{\frac{\partial {W}}{\partial t} +A(W)\frac{\partial W}{\partial x}=S(x,W)},
\label{cmp2mdl3}        
\end{align}
where 
\begin{align}
A(W)=
\begin{bmatrix}
U_{1}   & H_{1}     & 0               & 0  \\
g         &  U_{1}     & g               & 0   \\
0         &    0          & U_{2}        & H_{2}       \\
rg        &     0         & g              & U_{2}\end{bmatrix}
\label{matjac3}
\end{align}
Simple and exact analytical expression of the four eigenvalues of $A$ is
not obvious. Complicate expression 
can be obtained tediously by applying the Cardano-Vieta method.
For the case of $r\approx 1$ and $U_1\approx U_2$, a first order approximation of the eigenvalues is given in 
\cite{niecaspar11b, abgkar09}. 

In this work, we consider the case
where
$r\ll1$, namely, when  the bottom fluid is much thicker than the upper fluid.
Moreover, recall that our problem is to consider the boundary controller design to stabilize \eqref{cmp2mdl3}.
 \subsection{Linearization of the  \textit{Saint-Venant} model} 
 We denote the steady-state
associated to 
 the system  \eqref{cmp2mdl3} by $W^{*}=(H_{1}^{*},U_{1}^{*},H_{2}^{*},U_{2}^{*})$, which satisfies the following equation of a compact form:
 \begin{align}
 \bm{A(W^{*})\partial_x W^{*}=S(x,W^{*})}.
 \label{eqstead}
 \end{align}
 To obtain a constant steady-state, we work in the sequel with a flat bathymetry
$(\partial_x B=0)$.
A constant  steady-state of  the two-layer Saint-Venant
equations can be characterized by:
\begin{align}
\left\{
\begin{array}{llll}
 H_{1}^{*}U_{1}^{*}=cste,\qquad H_{2}^{*}U_{2}^{*}=cste,\\
 \frac{U_{1}^{*2}}{2}+g(H_{1}^{*}+H_{2}^{*})= -C_f|U_{1}^{*}-U_{2}^{*}|(U_{1}^{*}-U_{2}^{*})
, \\
\frac{U_{2}^{*2}}{2}+g(H_{2}^{*}+rH_{1}^{*})= rC_f|U_{1}^{*}-U_{2}^{*}|(U_{1}^{*}-U_{2}^{*}).
\end{array}
\label{steadystate2}\right.
\end{align}
In order to linearize the governing equations around the steady state, we define the deviation $(h_1,
u_1, h_2, u_2)$ of the state 
$(H_1, U_1, H_2, U_2)$ with respect to the steady-state  $(H_{1}^{*},U_{1}^{*},H_{2}^{*},U_{2}^{*})$
by:
\begin{align}
\left\{
\begin{array}{llll}
h_1=H_1-H_{1}^{*},\qquad  u_1=  U_1-U_{1}^{*},\\
h_2=H_2-H_{2}^{*},\qquad  u_2=  U_2-U_{2}^{*}.
\end{array}
\label{perturb1}\right.
\end{align}
Then, the linearized version of (\ref{cmp2mdl3})
can be written in  a matrix form as 
 \begin{align}
\mathbf{\partial_t U +A^*\partial_x U= S_{l}(U)},
\label{linversionmodel} 
\end{align}
 where
 $$U=(h_1, u_1, h_2, u_2)^{T}, \qquad  A^*=A(W^*),  $$
 and
 $$ S_{l}(U)=[0\; -\alpha_{s}^{f}(u_1-u_2)\; 0 \; r\alpha_{s}^{f}(u_1-u_2)]^T$$
with
  $$\alpha_{s}^{f}=2C_f|U_{1}^{*}-U_{2}^{*}|.$$
 We consider a constant steady state here for the sake of readability and simplicity in the presentation of the linear model.
\subsection{Linearized  \textit{Saint-Venant} model in Riemann coordinates}
We are  to explore the system eigenstructure of the linear form
(\ref{linversionmodel}) in this subsection. 
The characteristic equation derived from the matrix  $A^{*}$ is
\begin{align}
&\Theta=rg^2H_{1}^{*}H_{2}^{*},
\label{characpoly}
\end{align}
where
\begin{align}
&\Theta=\Big((\lambda-U_{1}^{*})^2-gH_{1}^{*}\Big)\Big((\lambda-U_{2}^{*})^2-gH_{2}^{*}\Big).
\label{characpoly}
\end{align}
For the case $r=0$, straightforward calculations lead to the following real
eigenvalues for $A^*$:
\begin{align}
\begin{array}{lll}
\lambda_{1}=U_{1}^{*}-\sqrt{gH_{1}^{*}}, \qquad
 \lambda_{2}=U_{1}^{*}+\sqrt{gH_{1}^{*}},\\
\lambda_{3}=U_{2}^{*}-\sqrt{gH_{2}^{*}}, \qquad
 \lambda_{4}=U_{2}^{*}+\sqrt{gH_{2}^{*}}.
\end{array}
\label{appvalpro}
\end{align}
We notice that the eigenvalues in this case are those corresponding to each
layer separately.
Following the results in \cite{jbjc53}, the eigenvalues for the system (\ref{cmp2mdl3})  in the 
case of $r\ll 1$ i.e $\rho_1\ll \rho_2$ approach to those given in (\ref{appvalpro}).
From (\ref{appvalpro}),  the internal and external characteristics  travel
at different speeds, and indeed, the 
lower layer characteristics moves much slower than the upper ones in
the case
of $r \ll 1$.
Let us now recast  the equation (\ref{linversionmodel})   into a diagonal form.
For a given eigenvalue $\lambda_k$ $(k=1,2,3,4)$ of the matrix $A^*$, the
associated right eigenvector is expressed  by
\begin{align}
V_k=
\begin{bmatrix}
1\\
\frac{\lambda_k-U_{1}^{*}}{H_{1}^{*}}\\
\frac{(\lambda_k-U_{1}^{*})^{2}-gH_{1}^{*}}{gH_{1}^{*}}\\
\frac{(\lambda_k-U_{2}^{*})((\lambda_k-U_{1}^{*})^{2}-gH_{1}^{*})}{gH_{1}^{*}H_{2}^{*}}
\end{bmatrix}
\label{rightvecpro}
\end{align}
 Some computations lead to the associated left eigenvector $L_k$, which is given
by:
\begin{align}
L_k=-\frac{1}{(\lambda_i-\lambda_k)(\lambda_j-\lambda_k)(\lambda_l-\lambda_k)}\nonumber\\
\times
\Big[l_{k,1}\quad l_{k,2}\quad l_{k,3}\quad l_{k,4}\Big]^{T}
\label{leftvectpro}
\end{align}
for  $i\neq j \neq l \neq k \in\{1,2,3,4\},$
where
\begin{align}
l_{k,1}&=U_{1}^{*3}-(trA^*-\lambda_{k})(U_{1}^{*2}+gH_{1}^{*})+f_k\nonumber\\
& +3gH_{1}^{*}-\frac{det A^*}{\lambda_k}, \\
l_{k,2}&=3H_{1}^{*}U_{1}^{*2}-2H_{1}^{*}U_{1}^{*}(trA^*-\lambda_{k})\nonumber\\
& + H_{1}^{*}(f_k+gH_{1}^{*}),\\
l_{k,3}&=gH_{1}^{*}(7U_{1}^{*}-\lambda_{k}), \quad
l_{k,4}=gH_{1}^{*}H_{2}^{*}.
\label{compleftvec}
\end{align}
The quantities $f_k$ are defined by:
  \begin{align}
  f_1=(\lambda_3+\lambda_2)\lambda_4 + \lambda_2\lambda_3,\,
  f_2=(\lambda_3+\lambda_1)\lambda_4 + \lambda_1\lambda_3,\\
  f_3=(\lambda_2+\lambda_1)\lambda_4 + \lambda_1\lambda_2,\,
  f_4=(\lambda_1+\lambda_2)\lambda_3 + \lambda_1\lambda_2.
  \end{align}
 We are to express the linear version  (\ref{linversionmodel}) of the governing equations in term of
the characteristic coordinates or
 Riemann Invariants. Multiplying the equation (\ref{linversionmodel}) by
the left eigenvectors 
 $L_k$  (each for a given eigenvalue $\lambda_k$) of the matrix $A^*$,  
 we get  that the characteristic coordinates (Riemann Invariants) $\xi_{k}$
 are:
\begin{align}
\xi_{k}=L^{tr}_{k}U
&=-\frac{1}{(\lambda_i-\lambda_k)(\lambda_j-\lambda_k)(\lambda_l-\lambda_k)},\nonumber\\
& \times \Big [l_{k,1}h_1+l_{k,2}u_1+l_{k,3}h_2+l_{k,4}u_2 \Big]
\end{align}
for $i\neq j \neq l \neq k \in \{1,2,3,4\}$. Therefore, we can express the variables $h_1$, $u_1$, $h_2$ and  $u_2$ in
term of the Riemann coordinates:
\begin{align}
\left\{
\begin{array}{llll}
h_1=\xi_{1}+\xi_{2}+\xi_{3}+\xi_{4},\\
u_1=\gamma_1\xi_{1}+\gamma_2\xi_{2}
+\gamma_3\xi_{3}+\gamma_4\xi_{4},\\
h_2=\beta_1\xi_{1}+\beta_2\xi_{2}+\beta_3\xi_{3}+\beta_4\xi_{4},\\
u_2=\alpha_1\xi_{1}+\alpha_2\xi_{2}+\alpha_3\xi_{3}+\alpha_4\xi_{4},
\end{array}
\label{inverserim}\right.
\end{align}
where
\begin{align}
&\gamma_k=\frac{\lambda_k-1}{H_{1}^{*}},\\
&\beta_k=\frac{1}{gH_{1}^{*}}\Big(U_{1}^{*2}+2(\lambda_k-1)U_{1}^{*}-\lambda_{k}^{2}+gH_{1}^{*}\Big),
\end{align}
and
\begin{align}
\alpha_k=&\frac{1}{gH_{1}^{*}H_{2}^{*}}\left((gH_{1}^{*}\beta_k-2\lambda_{k}^{2})U_{2}^{*}
+3U_{1}^{*3}\right. \nonumber\\
&\left.+7(\lambda_{k}-1)U_{1}^{*2}+2(gH_{1}^{*}-2\lambda_{k}^{2})U_{1}^{*}\right.\nonumber\\
&\left.
+\lambda_{k}^{2}(trA^*-\lambda_{k})+gH_{1}^{*}(\lambda_{k}+2)
\right).
\end{align}
We introduce the following more compact notations:
\begin{align}
\bm{\xi} =  \bpm \xi_1 &  \xi_2 &  \xi_3 &  \xi_4  \epm^{tr}, \;\;
\end{align}
and
\begin{align}
{\bm \Lambda} = \text{diag} \{ \lambda_1, \lambda_2, \lambda_3, \lambda_4
 \}.
\end{align}
Using the characteristic coordinates, we  recast the equation (\ref{linversionmodel})
into the following form:
 \begin{align}
\bm{\partial_t \xi +\Lambda\partial_x \xi=M \xi,}
\label{caractk0} 
\end{align}
where 
\begin{align}
&M(W^*)=\bpm 0 & \alpha_{s}^{f} & 0 &-r\alpha_{s}^{f} \epm^{tr}\nonumber\\
& ~~~~~~~~~~~~~~\times \bpm \alpha_1-\gamma_1 & \alpha_2-\gamma_2 & \alpha_3-\gamma_3 &
\alpha_4-\gamma_4 \epm.
\label{matrhs1}
\end{align}
 We consider the case where both layers have  a subcritical flow regime.
  Define the state vectors $$u(t,x)=(\xi_2, \xi_4), ~v(t,x)=(\xi_1, \xi_3),$$ and
introduce the transport speed matrices $${\bm \Lambda}^{\rm r} = \text{diag} \{
\lambda_2, \lambda_4  \},~{-\bm \Lambda}^{\rm l} = \text{diag} \{ \lambda_1, \lambda_3\}.$$
Then,  the system  (\ref{caractk0}) can be rewritten
as
 \begin{align}
\partial_t u(t,x) +\Lambda^{\rm r}\partial_x u(t,x) &= S^{\rm r}u(t,x)
\nonumber\\
&~~+ S^{\rm l}v(t,x),\label{caractknew-10}\\
\partial_t v(t,x) - \Lambda^{\rm l}\partial_x v(t,x) &= 0,
\label{caractknew-20}    
\end{align}
where
\begin{align}
S^{\rm r}&=
\begin{bmatrix}
\alpha_{s}^{f}(\alpha_1-\gamma_1)        &  \alpha_{s}^{f}(\alpha_2-\gamma_2)
 \\
r\alpha_{s}^{f}(\gamma_1-\alpha_1)        &  r\alpha_{s}^{f}(\gamma_2-\alpha_2)
 \end{bmatrix},\\ 
S^{\rm l}&=
\begin{bmatrix}
\alpha_{s}^{f}(\alpha_3-\gamma_3)        &  \alpha_{s}^{f}(\alpha_4-\gamma_4)
 \\
r\alpha_{s}^{f}(\gamma_3-\alpha_3)        &  r\alpha_{s}^{f}(\gamma_4-\alpha_4)
\end{bmatrix}
\label{matrhs20}
\end{align}
%
To close the writing of the system (\ref{caractknew-10})-(\ref{caractknew-20}),
we enclose to it the following boundary and initial conditions:
\begin{align}   
& u(t,0)=Q_0v(t,0)  \text{  and }  v(t,1)=R_1u(t,1) +\mathfrak{U}(t),\label{U0}\\
& u(0,x)=u_0(x) \text{  and }   v(0,x)=v_0(x),
 \label{bcindsyst10}
\end{align}
where $Q_0$,  $R_1$ $\in \mathcal{M}_{2,2}(\mathbb{R}),$ and $\mathfrak{U}(t)=(\mathfrak{u}_1(t),\mathfrak{u}_2(t))$ consists of the boundary controllers we need to design.
\section{controller design of a  general system}
In this section, we consider the backstepping controller design of a more general system, which could includes the \textit{Saint-Venant}  model as a special case. While solving our problem with the \textit{Saint-Venant}  model, it is also worth noting that the result derived in this section could be treated as a full theoretical result by itself.

\subsection{A more general control system}
The more general system discussed in this section is
 \begin{align}
&\partial_t u(t,x) +\Lambda^{\rm r}(x)\partial_x u(t,x)\nonumber\\
 &~~~~~~~~~~~~= S^{\rm r}(x)u(t,x)
+ S^{\rm l}(x)v(t,x),\label{caractknew-1}\\
&\partial_t v(t,x) - \Lambda^{\rm l}(x)\partial_x v(t,x) =S^{\rm o}(x)u(t,x),
\label{caractknew-2}    
\end{align}
where
        \begin{align}
                &u(x, t) = \begin{bmatrix} u_1(x, t), \, u_2(x, t), \,   \dots,
\, u_n(x, t) \\ \end{bmatrix}, \\
                &v(x, t) = \begin{bmatrix} v_1(x, t), \, v_2(x, t), \, \dots,
\,v_m(x, t) \\ \end{bmatrix}
        \end{align}
        are the systems states. The matrices       
         \begin{align}
                        &{\bm{\Lambda^{\rm r}}}(x)
                 = \text{diag} \begin{bmatrix} \lambda^{\rm r}_1(x)  \, ,\,
\lambda^{\rm r}_2(x)   \, ,\, \cdots  \, ,\,  \lambda^{\rm r}_{m-1}(x)  \end{bmatrix},\\
&{\bm{\Lambda^{\rm l}}} (x)
                 = \text{diag} \begin{bmatrix} \lambda^{\rm l}_1(x)  \, ,\,
\lambda^{\rm l}_2(x)   \, ,\, \cdots  \, ,\,  \lambda^{\rm l}_{m-1}(x)  \end{bmatrix},
        \end{align}
 subject to the restriction
        \begin{align}
                - \lambda^{\rm l}_m(x) <- \lambda^{\rm l}_2(x)\dots < - \lambda^{\rm
l}_1(x) < 0, \\
0<\lambda^{\rm r}_1(x) < \lambda^{\rm r}_2(x) < \dots <\lambda^{\rm r}_n(x),
        \end{align}and the in-domain parameters are given as
        \begin{align}
                        S^{\rm r}(x)
                & =
                                \{ S^{\rm r}_{ij}(x) \}_{1 \leq i \leq n,
1 \leq j \leq n} ,\\
S^{\rm l}(x)
                & =
                                \{ S^{\rm l}_{ij}(x) \}_{1 \leq i \leq n,
1 \leq j \leq m},\\
S^{\rm o}(x)
                & =
                               \{ S^{\rm o}_{ij}(x) \}_{1 \leq i \leq n,1 \leq j \leq m}. 
        \end{align}
 The system is also equipped with the following boundary and initial conditions:
\begin{align}   
& u(t,0)=Q_0v(t,0)  \text{  and }  v(t,1)=R_1u(t,1) +\mathfrak{U}(t),\label{U}\\
& u(0,x)=u_0(x) \text{  and }   v(0,x)=v_0(x),
 \label{bcindsyst1}
\end{align}
where the  boundary parameters $Q_0$,  $R_1$ $\in \mathcal{M}_{m,n}(\mathbb{R})$ are given as 
        \begin{equation}
                        R_1
                = 
                                \{ r_{ij} \}_{1 \leq i \leq m, 1 \leq j \leq
n}.
        \end{equation}
\begin{align}
                        Q_0
                & =
                                \{ q_{ij}(x) \}_{1 \leq i \leq n, 1 \leq
j \leq m,} 
        \end{align}  and the boundary controllers
are given as \begin{align}\mathfrak{U}(t)=\begin{bmatrix} \mathfrak{u}_1(t)  &\mathfrak{u}_2(t)& \dots & \mathfrak{u}_n(t)\end{bmatrix}^T.\end{align}

\subsection{Target system}
We employ the PDE backstepping method. First, we  construct a backstepping transformation to map  the system  (\ref{caractknew-1})-(\ref{caractknew-2}) into a  target system with  desirable stability property, which follows from the one constructed in   \cite{LongHu2015}.

   Consider the following target system
  \begin{align}
&\partial_t \epsilon (t,x) +\Lambda^{\rm r}(x)\partial_x \epsilon(t,x)
= S^{\rm r}(x)\epsilon(t,x) + S^{\rm l}(x)\beta(t,x)\nonumber\\
 &   +\int_{0}^{x}C^{\rm r}(x,\xi)\epsilon(\xi) d\xi  +\int_{0}^{x}C^{\rm l}(x,\xi)\beta(\xi) d\xi \label{targeteq0}      \\
&\partial_t \beta(t,x) - \Lambda^{\rm l}(x)\partial_x \beta(t,x) = \Delta(x)\beta(0,t)
\label{targeteq1}       
\end{align}

 with the following boundary conditions
  \begin{align} 
 \epsilon(t,0)=Q_0\beta(t,0)  \text{  and }  \beta(t,1)=0,
 \label{bcindtarget1}
\end{align}
where
 \begin{align}
\Delta(x)= \begin{bmatrix} 0 & \cdots & \cdots & 0\\ 
\delta_{2,1}(x) & \ddots & \ddots &\vdots\\
\vdots & \ddots & \ddots&\vdots\\
\delta_{m,1}(x)& \cdots & \delta_{m,m-1}(x) & 0  \end{bmatrix},\label{Delat_matr0}
\end{align} 
and $C^{\rm r}, ~C^{\rm l}$ are $L^{\infty}$ matrix functions  defined on the triangular domain
$$\mathbb{T}=\Big\{ (x,\xi)\in \mathbb{R}^2 | \; 0 \leq \xi  \leq x  \leq 1 \Big\}.$$
 
Here, $\Delta(x),~ C^{\rm r}, ~C^{\rm l}$ are all to be determined by introducing a backstepping transformation later.
\subsection{Backstepping controller design}
In order to map  the system  (\ref{caractknew-1})-(\ref{caractknew-2}) into the desired target system  (\ref{targeteq0})-(\ref{bcindtarget1}), we consider the following backstepping 
transformation
\begin{align}
\begin{pmatrix}
\epsilon(t,x)\\
\beta(t,x)
 \end{pmatrix}
=
\begin{pmatrix}
u(t,x)\\
v(t,x)
 \end{pmatrix}
-
\int_0^x \mathfrak{G}(x,\xi)
\begin{pmatrix}
u(t,\xi)\\ v(t,\xi)
 \end{pmatrix}
d\xi,
\label{backtransfo1}
\end{align}
where
\begin{align}
\mathfrak{G}=
\begin{pmatrix}
0 & 0\\
G_{21}(x,\xi)  & G_{22}(x,\xi)
 \end{pmatrix}.
\end{align}
Here the to-be-determined kernels $G_{21}$ and $G_{22}$ are defined on the domain $\mathbb{T}$, which
  need to satisfy the following system of equations:
\begin{align}\label{kerneleq1}
&\partial_\xi G_{21}(x,\xi){\bm{\Lambda^{\rm r}}}(\xi) - {\bm{\Lambda^{\rm l}}}(x)\partial_x G_{21}(x,\xi)\nonumber\\
&=\!\!-G_{21}(x,\xi) \frac{d{\bm{\Lambda^{\rm r}}}(\xi)}{d\xi}\!\!-\!\! G_{21}(x,\xi)S^{\rm r}(\xi)\!\!-\!\!G_{22}\!(x,\xi)S^{\rm o}\!(\xi)  \\
&\partial_\xi  G_{22}(x,\xi){\bm{\Lambda^{\rm r}}}(\xi)+ {\bm{\Lambda^{\rm l}}}(x)\partial_x  G_{22}(x,\xi)\nonumber\\
 &~~~~~=- G_{22}(x,\xi)\frac{d{\bm{\Lambda^{\rm r}}}(\xi)}{d\xi} + G_{21}(x,\xi)S^{\rm l}(\xi),\label{kerneleq2} 
\end{align}
and the following boundary conditions:
\begin{align}
G_{21}(x,x){\bm{\Lambda^{\rm r}}}(x)+{\bm{\Lambda^{\rm
l}}}(x) G_{21}(x,x)&=-S^{\rm o}(x), \label{kernelbc1}\\
 G_{22}(x,x){\bm{\Lambda^{\rm
l}}}(x) - {\bm{\Lambda^{\rm
l}}}(x) G_{22}(x,x)&=0, \label{kernelbc2}\\
 G_{21}(x,0){\bm{\Lambda^{\rm
r}}}(0) Q_0-G_{22}(x,0){\bm{\Lambda^{\rm
l}}}(0)&=-\Delta(x).\label{kernelbc3}
\end{align}
The existence, regularity and invertibility of the backstepping tranformation could be obtained similarly as  \cite{LongHu2015}, which is omitted here due to page limit consideration. Moreover, $\delta_{i,j}(x)$ for  $i=\overline{2,~m}$ 
$j=\overline{1,~ m-1}$ can be  obtained  from its  inverse transformation.
and the following equations are obtained for $C^{\rm r}(x,\xi)$,
$C^{\rm l}(x,\xi)$:
\begin{align}
C^{\rm r}(x,\xi)&=   S^{\rm l}(x) G_{21}(x,\xi)\nonumber\\
&~~~+  \int_{\xi}^{x} C^{\rm l}(x,\eta)
G_{21}(\xi,\eta) \,d\eta,  \label{cpluseq}\\
C^{\rm l}(x,\xi)&=  S^{\rm l}(x)G_{22}(x,\xi)\nonumber\\
&~~~~+\int_{\xi}^{x}C^{\rm l}(x,\eta)G_{22}(\xi,\eta)
\,d\eta.\label{cmoinseq}
\end{align}

 Hence, the control  law $\mathfrak{U}(t)$ can be obtained by plugging the transformation  (\ref{backtransfo1}) into
(\ref{U}). Readily, $\beta(t,1)=0$ implies that
\begin{align}
&\mathfrak{U}(t)= -R_1u(t,1)\nonumber\\
 &~~~~~~~+ \int_{0}^{1}\Big[ G_{21}(1,\xi)u(t,\xi)+ G_{22}(1,\xi)v(t,\xi)\Big]\,d\xi.\
\label{controldef}
\end{align}

\subsection{Stability of the Target system}
Assume that there exist constants $M>0, \bar q>0$, such that
\begin{align}
&\hspace{-3mm} \|C^r(x, \xi) \|, \|C^l(x, \xi) \|,\| S^{\rm r}(x)\|,  \|
S^{\rm l}(x)
\| \leq M,
\  \nonumber\\
&~~~~~~~~~~~~~~~~~~~~~~~~~~~~~~~~~~~~\forall x\in [0,1], \xi \in [0,x],\\
&\|Q^T_0Q_0\|<\bar{q,}
\end{align}
where $\|\cdot\|$ denotes the $2$-norm, and denote
\begin{align}
&\min\{\lambda_i^r(x), \lambda_i^l(x);x\in[0,1],i=\overline{1,n}\}=\underline \lambda,\\
&\max\{\lambda_i^r(x), \lambda_i^l(x);x\in[0,1],i=\overline{1,n}\}=\overline \lambda.
\end{align}We first prove exponential stability of the target system
(\ref{targeteq0})-(\ref{bcindtarget1}). The novelty, compared with other results in this area (i.e., \cite{LongHu2015}), lies in the newly proposed Lyapunov function. We employ a Lyapunov function that needs to be successively determined.

\begin{lemma}
For any given initial data $(\epsilon^0,  \ \beta^0)^T$\\=$(\epsilon(0, \cdot),  \ \beta(0, \cdot))^T \, $ $\in \, \left(\mathcal{L}^{2}([0,1])\right)^{n+m}$ and under the assumption that
$C^{\rm r}, ~C^{\rm l}  \, \in \mathcal{C}(\mathbb{T})$, the equilibrium $(\epsilon, \  \ \beta)^T=(0, \ 0, \ 0, \ 0)^T$ of the target system
(\ref{targeteq0})-(\ref{Delat_matr0}) is  exponentially stable in the $\mathcal{L}^2$-norm:
 \begin{align}
&\|( \epsilon(t,\cdot),  \beta(t,\cdot))\|^2_{\mathcal{L}^2}\nonumber\\
&~~~~~~~:= \int_0^1 \epsilon^T(t,x)
\epsilon(t,x)+\beta^T(t,x) \beta(t,x)dx.
\end{align}

\end{lemma}

\textbf{PROOF:}\\
 We consider the  following Lyapunov function:
  \begin{align}
&V(t)=\frac{1}{2} \int_0^1 e^{-\nu x}\epsilon^T(t,x)\bm{\Lambda^{\rm r}_{\rm inv}}(x) \epsilon(t,x)dx\nonumber\\
&~~~~~~+\frac{1}{2} \int_0^1 (1+x)\beta^T(t,x) D\bm{\Lambda^{\rm l}_{\rm inv}}(x) \beta(t,x)  dx, \label{V1}
\end{align}
where
$
D= \text{diag} \begin{bmatrix}d_1 \, ,\,d_2  \, ,\, \cdots  \, ,\, d_{m-1}
 \, ,\, d_m \end{bmatrix}$, and
\begin{align*}
&{\bm \Lambda^{\rm r}_{\rm inv}(x)}  = \text{diag} \left\{ \frac{1}{\lambda^{\rm
r}_1(x)},\dots, \frac{1}{\lambda^{\rm r}_n(x)} \right\},\\
& {\bm \Lambda^{\rm l}_{\rm inv}(x)}= \text{diag} \left \{ \frac{1}{\lambda^{\rm
l}_1(x)},\dots, \frac{1}{\lambda^{\rm l}_m(x)}\right\}  
\end{align*}
 with $\left(\lambda^{\rm r}_1(x),\dots,\lambda^{\rm r}_n(x)\right)>0,~\left(\lambda^{\rm l}_1(x),\dots, \lambda^{\rm r}_m(x)\right)>0. $ The constants  $\nu$ and $d_{1},\, d_2  \, ,\, \cdots  \, ,\, d_{m-1}$ are all positive parameters 
to be determined. Then, we have
\begin{align}
C_1\!\|( \epsilon(t,\cdot),  \beta(t,\cdot))\|^2_{\mathcal{L}^2}\!\!\leq V(t)\leq\!\! C_2 \|( \epsilon(t,\cdot),  \beta(t,\cdot))\|^2_{\mathcal{L}^2},\label{equivalent}
\end{align}
where the two positive constants are
\begin{align}
&C_1=\frac{1}{2\bar \lambda} \min\left\{e^{-\nu}, ~ d_{i}; i\,
=\overline{1,m} \right\},\label{c1}\\
&C_2=\frac{1}{2\underline \lambda} \max\left\{1, ~  2d_{i}; i\,
=\overline{1,m}\right\}.\label{c2}
\end{align}
This ensures that $V(t)$ is positive definite.
 
Differentiating \eqref{V1} with respect to time, we get:
\begin{align}
& \dot V(t)= \int_0^1 e^{-\nu x} \epsilon^T(t,x) \bm{\Lambda^{\rm r}_{\rm inv}}(x)  \partial_t  \epsilon(t,x) dx\nonumber\\
& ~~~~~~~+ \int_0^1 (1+x)\beta^T(t,x)D\bm{\Lambda^{\rm l}_{\rm inv}}(x) \partial_t  \beta(t,x)dx \hspace{-0.4mm}.\label{dotV1-00}
\end{align}
With the help of (\ref{targeteq0}) and (\ref{targeteq1}),  we rewrite  \eqref{dotV1-00}  as follows:
\begin{align}
\hspace{-0.07in} \dot V(t)=& \int_0^1 -e^{-\nu x} \epsilon^T(t,x)\partial_x \epsilon(t,x)dx \nonumber\\ &+\int_0^1 e^{-\nu x} \epsilon^T(t,x)\bm{\Lambda^{\rm r}_{\rm inv}}(x) S^{\rm r}(x)\epsilon(t,x) dx \nonumber\\ 
&+\int_0^1 e^{-\nu x} \epsilon^T(t,x)\bm{\Lambda^{\rm r}_{\rm inv}}(x) S^{\rm l}(x)\beta(t,x)dx \nonumber\\
 &+ \int_0^1  e^{-\nu x}\epsilon^T(t,x) \left(\int_{0}^{x} \bm{\Lambda^{\rm r}_{\rm inv}}(x) C^{\rm r}(x,\xi)\epsilon(\xi) d\xi\right)dx\nonumber\\
 &+\int_0^1 e^{-\nu x}\epsilon^T(t,x) \left(\int_{0}^{x}\bm{\Lambda^{\rm r}_{\rm inv}}(x)C^{\rm l}(x,\xi)\beta(\xi) d\xi \right)dx \nonumber  \\
&+\int_0^1  (1+x)\beta^T(t,x)D\partial_x \beta(t,x)dx+ I(t),\label{dotV1-1}
\end{align}
where
\begin{align}
I(t)=\int_0^1 (1+x) \beta^T(t,x) D\bm{\Lambda^{\rm l}_{\rm inv}} (x) \Delta(x)\beta(0,t)\,dx.
\end{align}

Since
\begin{align*}
&I(t) \leq
(m-1)\int_0^1 \frac{(1+x)}{2}\beta^T(t,x)\beta(t,x) \, dx\nonumber\\ &+\beta_1(t,0)^2\int_0^1 \frac{(1+x)}{2} \sum_{i=2}^{m}d_i^2 \frac{1}{\lambda^{\rm
l}_i(x)^2}\delta_{i,1}^2(x)dx\nonumber\\
&+\beta_2(t,0)^2\int_0^1 \frac{(1+x)}{2} \sum_{i=3}^{m}d_i^2
\frac{1}{\lambda^{\rm
l}_i(x)^2}\delta_{i,2}^2(x)dx\nonumber\\ 
&+\cdots\nonumber\\ 
&+\beta_{m-2}(t,0)^2\int_0^1 \frac{(1+x)}{2}\sum_{i=m-1}^{m}d_i^2
\frac{1}{\lambda^{\rm
l}_i(x)^2}\delta_{i,m-2}^2(x)dx\nonumber\\
& +\beta_{m-1}(t,0)^2\int_0^1 \frac{(1+x)}{2}d_m^2
\frac{1}{\lambda^{\rm
l}_m(x)^2}\delta_{m,m-1}^2(x)dx,
\end{align*}
then
\begin{align}
&\dot V(t)\leq J_1(t)+J_2(t),
\end{align}
where
\begin{align}
&J_1(t)\nonumber\\
&=\beta_1(t,0)^2\nonumber\\
&~~~\times \Bigg \{\frac{\bar q}{2} - \frac{1}{2}d_1+\int_0^1
\frac{(1+x)}{2} \sum_{i=2}^{m}d_i^2 \frac{1}{\lambda^{\rm
l}_i(x)^2}\delta_{i,1}^2(x)dx\Bigg \} \nonumber\\
&+\beta_2(t,0)^2\nonumber\\
&\times \Bigg \{\frac{ \bar q}{2}-   \frac{1}{2}d_2+\int_0^1 \frac{(1+x)}{2} \sum_{i=3}^{m}d_i^2
\frac{1}{\lambda^{\rm
l}_i(x)^2}\delta_{i,2}^2(x)dx\Bigg \}\nonumber\\
&+\cdots+\beta_{m-2}(t,0)^2 \Bigg \{\frac{
\bar q}{2}-   \frac{1}{2}d_{m-2}\nonumber\\
&~~~~~~~~~~~~~~~+\int_0^1 \frac{(1+x)}{2}\sum_{i=m-1}^{m}d_i^2
\frac{1}{\lambda^{\rm
l}_i(x)^2}\delta_{i,m-2}^2(x)dx \bigg\}\nonumber\\
& +\beta_{m-1}(t,0)^2\Bigg \{\frac{ \bar q}{2}-   \frac{1}{2}d_{m-1} \nonumber\\ &~~~~~~~~~~~~~~~~~~+ \int_0^1 \frac{(1+x)}{2}d_m^2
\frac{1}{\lambda^{\rm
l}_m(x)^2}\delta_{m,m-1}^2(x)dx \Bigg\}\nonumber\\
& +\beta^2_m(t,0)\bigg \{\frac{ \bar q}{2}-   \frac{1}{2}d_{m}\bigg\},  \label{dotV1-1}
\end{align}
and 
\begin{align}
&J_2(t)=-\frac{1}{2}f_1(\nu)\int_0^1 e^{-\nu x} \epsilon^T(t,x) \epsilon(t,x)dx\nonumber\\
&~~~~~~~~~~~-\frac{1}{2}f_2(\nu)\int_0^1 \beta^T(t,x)\beta(t,x) \, dx
\end{align}
with
\begin{align}
&f_1(\nu)=\nu-2 \left(\frac{M}{\underline
\lambda}\right)^2-\frac{M}{\underline
\lambda}(5+\frac{1}{\nu}),\\
&f_2(\nu)=\min\left\{d_{i}; i\, =\overline{1,m}\right\}-2m+1-\frac{1}{\nu}.
\end{align}
First, choose the positive constants $d_{1},\, d_2  \, ,\,
\cdots  \, ,\, d_{m-1}$ successively as follows:
\begin{align}
&d_{m}\geq \bar q,\nonumber\\
&d_{m-1}\geq  \bar q+\int_0^1(1+x)d_m^2
\frac{1}{\lambda^{\rm
l}_m(x)^2}\delta_{m,m-1}^2(x)dx,\nonumber\\
&d_{m-2}\geq  \bar q+\int_0^1(1+x)\sum_{i=m-1}^{m}d_i^2
\frac{1}{\lambda^{\rm
l}_i(x)^2}\delta_{i,m-2}^2(x)dx,\nonumber\\
&~~~~~~~~~~~~~~~~~~~~~~\vdots\nonumber\\
&d_2\geq  \bar q+ \int_0^1(1+x) \sum_{i=3}^{m}d_i^2
\frac{1}{\lambda^{\rm
l}_i(x)^2}\delta_{i,2}^2(x)dx,\nonumber\\
&d_1\geq  \bar q+ \int_0^1(1+x) \sum_{i=2}^{m}d_i^2 \frac{1}{\lambda^{\rm
l}_i(x)^2}\delta_{i,1}^2(x)
\, dx,
\end{align}
which guarantee that $J_1(t)$ is non-negative.  Then, choose $\nu>0$ large enough to satisfy
\begin{align}
&f_1(\nu)>0,\\
&f_3(\nu):=\bar q-2m+1-\frac{1}{\nu}>0.
\end{align}\\ 
%
 from which we have
\begin{align}
&f_2(\nu)\geq f_3(\nu)>0. 
\end{align}
Thus, it holds that
\begin{align}
&\dot V(t)\leq J_2(t)\leq -cV(t),
\end{align}
where
\begin{align}
&c=\underline \lambda \min\left\{ f_1(\nu),~ \frac{1}{\max\{d_{i}; i\, =\overline{1,m}\}} f_2(\nu)   \right\},
\end{align}
which gives 
\begin{align}
V(t)\leq V(0)e^{-ct}.
\end{align}
Finally, from \eqref{equivalent}, we derive
\begin{align}
\|( \epsilon(t,\cdot),  \beta(t,\cdot))\|_{\mathcal{L}^2}
&\leq \sqrt{\frac{C_2}{C_1}}\|( \epsilon^0(\cdot),  \beta^0(\cdot))\|_{\mathcal{L}^2}e^{-ct},
\end{align} 
where $C_1, C_2$ are defined in \eqref{c1}, \eqref{c2}. This completes the proof.

\subsection{Stability of the closed-loop control system}
With the exponential stability of the target system, and with existence, regularity and invertibility of the backstepping tranformation, we are now ready to derive the stability of the closed-loop control system.
\begin{thm}\label{theo}
For any given initial data
$(u^0,  \ v^0)^T= \,$ $(u(0, \cdot),
 \ v(0, \cdot))^T\in \, \left(\mathcal{L}^{2}([0,1])\right)^{n+m}$
and under the assumption that
$C^{\rm r}, ~C^{\rm l}  \, \in \mathcal{C}(\mathbb{T})$, the equilibrium
$(u, \  \ v)^T=(0, \ 0, \ 0, \ 0)^T$ of the closed-loop system
(\ref{targeteq0})-(\ref{bcindtarget1}) with the designed controller (\ref{controldef})
is  $\mathcal{L}^2$-exponentially
stable.
\end{thm}
The exponential stability result stated in Theorem \ref{theo}  is immediatly applied to the linearized bilayer Saint-Venant  model \eqref{model30}, in Riemanian  coordinates, which consists  of  a coupled system of two rightward
and two
leftward convecting transport PDEs.
\section{Simulation results}
\begin{figure}[H]
\centering
\includegraphics[width=0.3\textwidth]{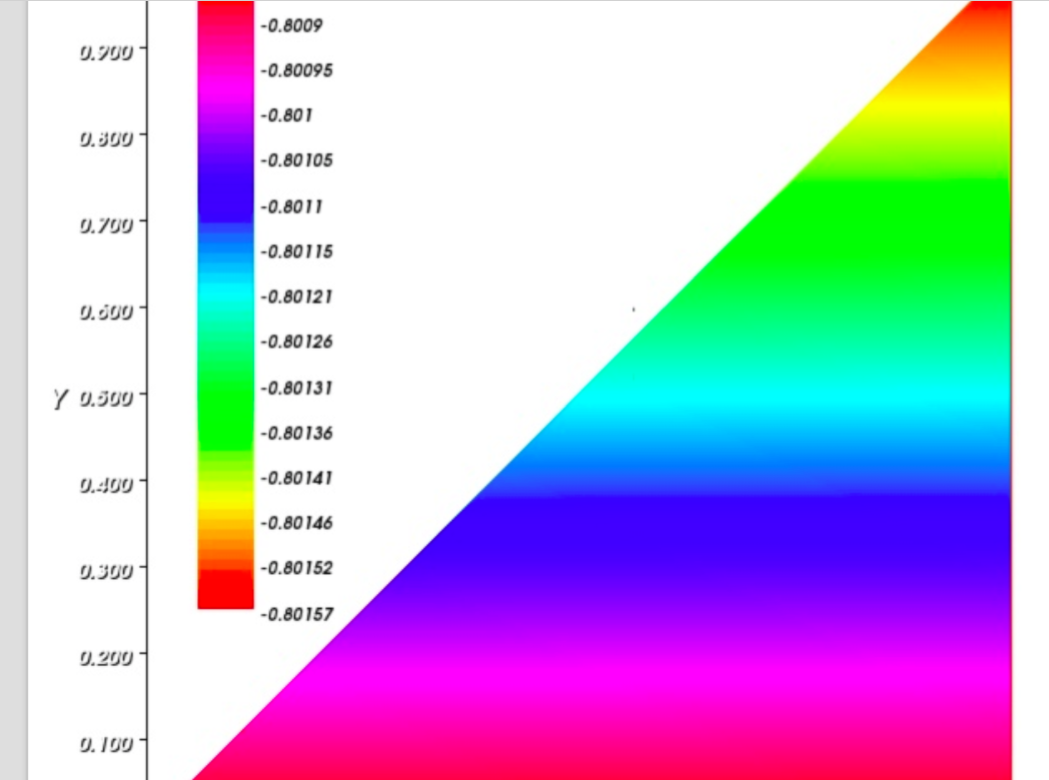}
\caption{The kernel $G_{21}$.}
\label{fig:solutionkernelG21}
\end{figure}

\begin{figure}[H]
\centering
\includegraphics[width=0.35\textwidth]{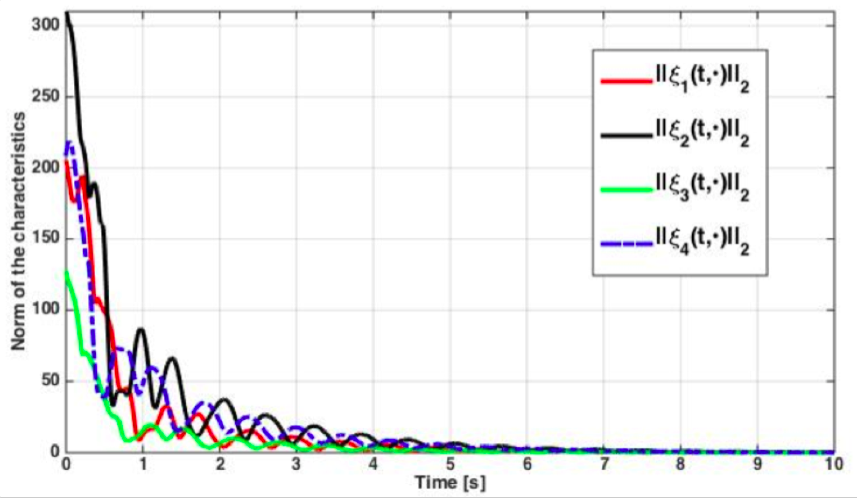}
\caption{Evolution in time of the norm of the characteristic solutions.}
\label{fig:normsolution-state}
\end{figure}

\begin{figure}[H]
  \centering
\includegraphics[width=0.35\textwidth]{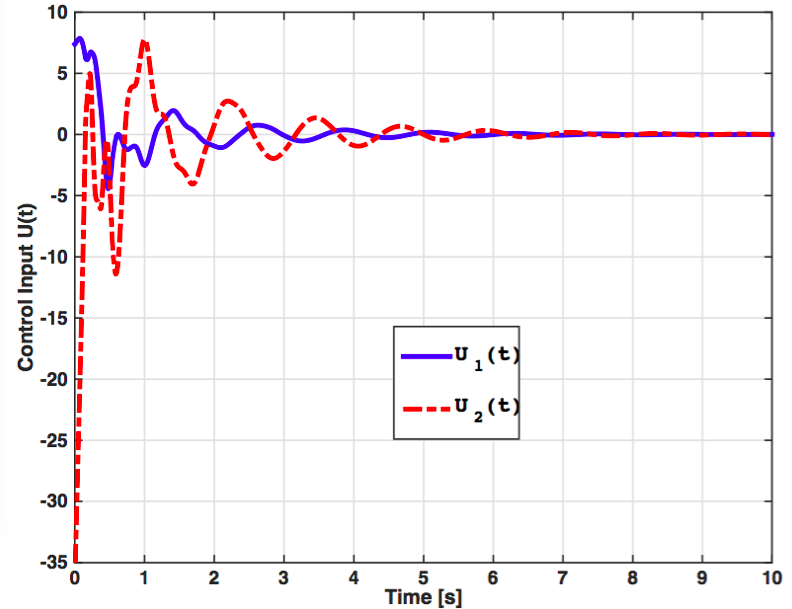}
      \caption{Evolution   of the component of the control input $U(t)$.}
         \label{fig:control-input-state}
\end{figure}

The goal of the following numerical simulations is to illustrate the efficiency of the designed
$\mathfrak{U}(t)$, namely (\ref{controldef}), to stabilize the linear system  (\ref{caractk0})
around the zero equilibrium. 
As initial conditions, the following data are considered for  the layer $1$ and $2$ through the physical variables
$$H_2(0,x)=2+0.5\exp\Big(-\frac{(x-0.5)^2}{0.003}\Big),\,H_1(0,x)=6-H_1(x)$$
and 
$$U_1(0,x)=\frac{10}{H_1(0,x)}+ 3\sin(2\pi x),$$ 
$$U_2(0,x)=-\frac{10}{H_2(0,x)}-3\sin(2\pi x).$$ 
The  initial data of the characteristic variables $\xi_k$, $(k=1, \,2, \,3, \,4)$ (for system (\ref{caractk0}))
 are computed as  function of the physical variables  $H_i(0,x)$ and   $U_i(0,x)$ for $i=1,\,2$,
thanks to the relation (\ref{leftvectpro}).
For the sake of simplicity, we consider the  following uniform steady state:
$$H^{*}_{1}=3,\; U^{*}_{1}=1,\, H^{*}_{2}=1,\, U^{*}_{2}=0.95.$$
With this choice of steady state  (set point), the  characteristic speeds are given by:
$$ \lambda_1 = 6.42,\, \lambda_2= 4.08,\,  \lambda_3 = -4.42  \text{ and }   \lambda_4 = -2.18.$$
Elsewhere, in the reported numerical experiments, the ratio $r$ between the densities is set to 
$0.01$ and the friction coefficient $C_f$ to $0.05$. We compute the solution up to time $T=10$.
Regarding to the boundary conditions (\ref{bcindsyst1}) the following matrix are considered
\begin{align}
Q_0 = \begin{bmatrix}
              -1.5 & 0.01 \\
                  0.01 & 1.5  \\
                \end{bmatrix}, \quad
R_1=\begin{bmatrix}
           0.5 & 0.1 \\
          0.15 & -0.5    \\
         \end{bmatrix}
        \end{align}
Our implementation is based on an accurate finite volume method for the evolutionary equation (\ref{caractk0}). 
More precisely, we use a modified Roe's scheme (see \cite{levequebook}). Since the computation of the designed 
control $\mathfrak{U}(t)$ requires the knowledge of the kernels $G_{21}$ and  $G_{22}$, the kernel problems 
are solved numerically according to \eqref{kerneleq1}-\eqref{kernelbc3}  using the finite element setup. 
As an illustration, the numerical solution of the second  component  of the kernel $G_{21}$ is given in Figure 
\ref{fig:solutionkernelG21}.
 Figure \ref{fig:normsolution-state} depicts  the evolution  in time of the $\mathcal{L}^2$-norm of the characteristics. 
As expected from  the theoretical part we observe clearly that the norm of all characteristics  decreases in time and 
converges to zero. As a result, this shows that the system  (\ref{caractk0}) converges to the zero equilibrium. 
Thereby the two-layer Saint-Venant model  (\ref{model30}) also converges to  ($H^{*}_{1}$, $U^{*}_{1}$, $H^{*}_{2}$, $U^{*}_{2}$).
In Figure \ref{fig:control-input-state} are depicted the behavior in time of each component of  the input control $\mathfrak{U}(t)$. 
Clearly, despite the initial amplitude of $\mathfrak{u}_2(t)$, this latter one decreases in time and vanishes after $t\geq4\, s$. 
Likewise  $\mathfrak{u}_2(t)$, the first component of the control input $\mathfrak{u}_1(t)$ shows the same trend with 
its amplitude decreasing in time and tending to zero after $t\geq 7s$ as can be seen in  Figure \ref{fig:control-input-state}.
As can be seen from these numerical simulations, the system (\ref{caractk0})  subject to the feedback control 
$\mathfrak{U}(t)$ is stabilized around the zero equilibrium as expected from the theoretical part.
 
\section{Conclusion and Future Works}
In this paper, a general system with spatially varying coefficients, consisting of  $m$ rightward and $n$ leftward convecting
transport PDEs, is first exponentially stabilized by $m$ backstepping boundary controllers.
Our   controller design idea is similar as an existing result for this system  with constant coefficients, and we employ a novel Lyapunov
function in which the parameters need to be successively determined in the stability proofs.
Then, this result is immediately applied to exponentially stabilize a 1D linearized bilayer
 \textit{Saint-Venant} model, which is
a coupled system of two rightward
and two
leftward convecting transport PDEs. One of our next steps is to consider the output feedback stabilizing of the general system, and then, as this paper, apply the result to the 1D bilayer
 \textit{Saint-Venant} problem.




\begin{thebibliography}{10}
\providecommand{\url}[1]{#1}
\csname url@samestyle\endcsname
\providecommand{\newblock}{\relax}
\providecommand{\bibinfo}[2]{#2}
\providecommand{\BIBentrySTDinterwordspacing}{\spaceskip=0pt\relax}
\providecommand{\BIBentryALTinterwordstretchfactor}{4}
\providecommand{\BIBentryALTinterwordspacing}{\spaceskip=\fontdimen2\font plus
\BIBentryALTinterwordstretchfactor\fontdimen3\font minus
  \fontdimen4\font\relax}
\providecommand{\BIBforeignlanguage}[2]{{%
\expandafter\ifx\csname l@#1\endcsname\relax
\typeout{** WARNING: IEEEtran.bst: No hyphenation pattern has been}%
\typeout{** loaded for the language `#1'. Using the pattern for}%
\typeout{** the default language instead.}%
\else
\language=\csname l@#1\endcsname
\fi
#2}}
\providecommand{\BIBdecl}{\relax}
\BIBdecl

\bibitem{Audusse2011}
\BIBentryALTinterwordspacing
E.~Audusse, M.-O. Bristeau, B.~Perthame, and J.~Sainte-Marie,
  ``\BIBforeignlanguage{eng}{A multilayer saint-venant system with mass
  exchanges for shallow water flows. derivation and numerical validation},''
  \emph{\BIBforeignlanguage{eng}{ESAIM: Mathematical Modelling and Numerical
  Analysis}}, vol.~45, no.~1, pp. 169--200, 1 2011. [Online]. Available:
  \url{http://eudml.org/doc/197421}
\BIBentrySTDinterwordspacing

\bibitem{bouchut2008}
F.~Bouchut and d.~L.~T. Morales, ``An entropy satisfying scheme for two-layer
  shallow water equations with uncoupled treatment,'' \emph{ESAIM: Mathematical
  Modelling and Numerical Analysis}, vol.~42, pp. 683--698, 7 2008.

\bibitem{Narbona2011}
G.~Narbona-Reina and J.~D.~D. Zabsonre, ``Existence of a global weak solution
  for a 2d viscous bi-layer shallow water model,'' \emph{Nonlinear Analysis:
  Real World Applications}, no.~10, pp. 2971--2984, 2009.

\bibitem{Munoz2003}
M.~L. Munoz-Ruiz, F.~J. Chatelon, and P.~Orenga,
  ``\BIBforeignlanguage{English}{On a bi-layer shallow-water problem},''
  \emph{\BIBforeignlanguage{English}{Nonlinear Analysis: Real World
  Applications}}, vol.~4, no.~1, pp. 139--171, 2003.

\bibitem{Dimeglio2013}
F.~Di~Meglio, R.~Vazquez, and M.~Krstic, ``Stabilization of a system of coupled
  first-order hyperbolic linear {PDE}s with a single boundary input,''
  \emph{IEEE Transactions on Automatic Control}, vol.~58, no.~12, pp.
  3097--3111, 2013.

\bibitem{krstic2008boundary}
M.~Krstic and A.~Smyshlyaev, \emph{Boundary control of PDEs: A course on
  backstepping designs}.\hskip 1em plus 0.5em minus 0.4em\relax Siam, 2008,
  vol.~16.

\bibitem{LongHu2015}
L.~Hu, F.~D. Meglio, R.~Vazquez, and M.~Krstic, ``Control of homodirectional
  and general heterodirectional linear cou- pled hyperbolic pdes,'' \emph{arXiv
  preprint arXiv:1504.07491}, 2015.

\bibitem{niecaspar11b}
E.~F. Nieto, M.~J. Castro-Diaz, and C.~Par\'es, ``On an intermediate field
  capturing riemann solver based on a parabolic viscosity matrix for the
  two-layer shallow water system,'' \emph{Journal of Scientific Computing},
  vol. Volume 48, Issue 1-3, pp. 117--140, July 2011.

\bibitem{abgkar09}
R.~Abgrall and S.~Karni, ``Two-layer shallow water system: a relaxation
  approach,'' \emph{SIAM J. Sci. Comput.}, vol. 31, No 3, pp. 1603--1627, 2009.

\bibitem{jbjc53}
J.~Schijf and J.~Schonfeld, ``Theoretical considerations on the motion of salt
  and fresh water,'' \emph{Proc. of the Minn. Int. Hydraulics Conv. Joint
  meeting IAHR and Hyd. Div. ASCE.}, pp. 321--333, Sept. (1953).

\bibitem{levequebook}
\BIBentryALTinterwordspacing
R.~J. LeVeque, \emph{Finite volume methods for hyperbolic problems}, ser.
  Cambridge texts in applied mathematics. Cambridge, New York: Cambridge University Press, 2002.
\BIBentrySTDinterwordspacing

\end{thebibliography}
\end{document}